\documentclass[10pt,a4paper]{article}
\usepackage[latin9]{luainputenc}
\usepackage{float}
\usepackage{booktabs}
\usepackage{amsmath}
\usepackage{amssymb}
\usepackage{graphicx}

\makeatletter

\@ifundefined{pageheight}{\let\pageheight\pdfpageheight}{}
\@ifundefined{pagewidth}{\let\pagewidth\pdfpagewidth}{}
\pageheight\paperheight
\pagewidth\paperwidth

\usepackage{amsfonts}
\usepackage{amsthm}
\usepackage{xcolor}
\usepackage{pbox}
\usepackage{setspace}
\usepackage{enumerate}
\newtheorem{theorem}{Theorem}[section]
\newtheorem{definition}[theorem]{Definition}\newtheorem{examples}[theorem]{Examples}\newtheorem{remark}[theorem]{Remark}\usepackage{qtree}

\makeatother

\begin{document}
	\begin{center}
		{\Large\bf{Lie point symmetries of the biharmonic heat equation on surfaces of revolution}}
	\end{center}
	\begin{center}
		{Aminu M. Nass$^{1,*}$, Kassimu Mpungu$^{2,a}$, Rahmatullah Ibrahim Nuruddeen$^{3,b}$  \\
		\medskip 
			\normalsize{$^{1}$Department of Actuarial Science, Federal University Dutse, Dutse,  P.M.B. 7156 Jigawa State, Nigeria.} \\
			\normalsize{$^{2}$Department of Physical Sciences, Mountains of the Moon University, Fort Portal, Uganda.} 
		   \normalsize{$^{3}$Department of Mathematics, Faculty of Physical Sciences, Federal University Dutse, Dutse, \\ P.M.B. 7156 Jigawa State, Nigeria.} \\
												\medskip 
			{Email addresses:  ${^*}$aminu.nass@fud.edu.ng, ${^a}$kmpungu@mmu.ac.ug, ${^b}$rahmatullah.n@fud.edu.ng } \\
			{Email addresses: ${^*}$Corresponding Author} } 
	\end{center}

\begin{abstract}
This paper uses Lie symmetry analysis to investigate the biharmonic heat equation on a generalized surface of revolution. We classify the Lie point symmetries associated with this equation, allowing for the identification of surfaces and the corresponding infinitesimal generators. In a significant move, we demonstrate that the biharmonic heat equation on a surface of revolution admits the same Lie symmetries as the harmonic heat equation on the same surface, highlighting a profound structural relationship between the two equations. Utilizing these symmetry groups, we derive similarity reductions that yield invariant forms of the equation and facilitate the construction of exact solutions. Finally, we provide certain examples illustrating precise solutions on the related surfaces with positive, negative, and zero Gaussian curvatures, demonstrating the versatility of the approach. This work contributes to the understanding of biharmonic heat equations on symmetric surfaces.
\\
\hspace{0.35 in}\bf{AMS Subject Classifications:}\rm \; 31A30; 35B06; 58J70.
 \end{abstract}

\section{Introduction}
The biharmonic equation is a fundamental partial differential equation with wide-ranging applications across scientific and engineering domains. It describes diverse phenomena, such as the bending of plates \cite{timoshenko}, fluid mechanics \cite{Batchelor}, heat conduction \cite{Ozisik}, image processing \cite{Bertalmio}, fluid film lubrication \cite{Hamrock}, electrostatics \cite{Griffiths}, geophysics \cite{Turcotte}, and the shape optimization \cite{Delfour} to mention but a few. Surfaces of revolution, due to their rotational symmetry, are widely used in engineering, industrial design \cite{mott}, and computer graphics \cite{foley}; thereby making them significant subjects for mathematical modeling and analysis. As for the golden methodology, Ovsiannikov \cite{ovsiannikov} transformed classical Lie symmetry theory into a more systematic framework, enabling its broad application to differential equations. In this context, Azad and Mustafa \cite{tahir}, addressed symmetry classification for the wave equation on a sphere and derived exact solutions via similarity reductions. Mpungu \cite{mpungu} extended this research by investigating the heat and wave equations on surfaces of revolution, deriving minimal symmetry algebras and classifying all surfaces with higher symmetry dimensions.  Additionally, the study encompassed symmetry reductions and corresponding exact solutions.  For other insightful studies on the classical Lie symmetry theory of differential equations, its applications, and extensions, we refer the reader to \cite{Azad, mpungu1, nass1, nass2, aminu, ibragimov, bluman, R1, R2, R3} and the references provided therein. Biharmonic curves have also been studied in Euclidean and Minkowski  space \cite{abbena, inoguchi}, as well as surfaces with constant curvature and surfaces of revolution \cite{piu}. Besides, while a general classification remains difficult, several special cases have been resolved. 

However, motivated by the aforementioned developments, this paper addresses the symmetry classification problem for the biharmonic heat equation on a surface of revolution. Using Lie group analysis, we identify the infinitesimal generators and classify symmetries of the equation. We delve into the concept of similarity reductions for each symmetry class. By reducing the number of independent variables in the equation, we pave the way for the derivation of exact invariant solutions. These exact solutions hold paramount significance as they reveal the diverse patterns and behaviors exhibited by the biharmonic equation on surfaces of revolution.
Moreover, the arrangement of the manuscript follows the following manner: in section 2, we present some necessary mathematical backgrounds. Section 3 discusses the Lie symmetries of the biharmonic heat equation on the surfaces of revolution, including their classifications. In section 4, we explore symmetry reductions and construct exact solutions for some selected surfaces. Finally, we present some examples and conclude with insights on the structural implications of the overall findings in Section 5.

\section{Mathematical Background}\label{s2}
We consider a unit speed curve parameterized as follows \cite{mpungu},
\begin{equation}\label{unit}
 \varphi(x) = (v(x), w(x)), \quad where \	v'(x)^2 + w'(x)^2 = 1.
\end{equation}
\begin{theorem}
Let $\alpha(x)=(v(x),w(x))$ be a unit speed curve with a regular parameterization. The surface of revolution generated by rotating this curve about $z-$axis is given by
	\begin{equation*}
	X=(v(x),w(x)\cos(y),w(x)\sin(y)),
	\end{equation*}
while the Gaussian curvature is expressed as follows 
	\begin{equation}\label{2.1}
		\mathcal{K}=-\frac{w{''}(x)}{w{'}(x)}.
	\end{equation}
\end{theorem}
Next, we focus on the surfaces generated by revolving a unit speed profile curve 
$$\alpha(x)=(v(x),w(x))=(v(x),e^{f(x)}),$$
parameterized by 
\begin{equation}
	X(x,y)=(v(x),e^{f(x)}\cos(y),e^{f(x)}\sin(y)), \qquad {0}\leqslant{y}\leqslant{2\pi},
\end{equation}
with the metric given by 
\begin{equation}\label{key}
	g=dx^2+e^{2f(x)}dy^2.
\end{equation}

\begin{definition}
Let $(S,g)$ be Riemannian manifolds with metric $g.$ Then, the Laplacian is defined as follows \cite{mpungu},
\begin{equation}\label{1r} 
\Delta u=\frac{1}{\sqrt{\mid g\mid }}\frac{\partial }{\partial {{x}^{i}}}\left( \sqrt{\mid g\mid }{{g}^{ij}}\frac{\partial }{\partial {{x}^{j}}}(u) \right).
\end{equation}
\end{definition} 

\begin{definition}
The governing biharmonic heat equation on a surface of revolution through the latter Laplacian \eqref{1r} and the  metric (\ref{key}) takes the following expression
\begin{equation}\label{key1}
u_t=\left(f^{'}(x)\frac{\partial}{\partial x}+\frac{\partial^2}{\partial x^2}+e^{-2f(x)}\frac{\partial^2}{\partial y^2}\right)^{2}u.
\end{equation} 
\end{definition}

\section{Lie Symmetry Method}
To analyze the biharmonic heat equation on a surface of revolution \eqref{key1}, we employ the classical Lie group method. We assume the equation is invariant under the following one-parameter $(\epsilon)$ group transformations 
\begin{equation}\label{one}
\begin{split}
& \bar{x}=x+\epsilon\xi(x,y,t,u)+O(\epsilon^2), \\
& \bar{y}=y+\epsilon\phi(x,y,t,u)+O(\epsilon^2), \\
& \bar{t}=t+\epsilon\tau(x,y,t,u)+O(\epsilon^2), \\
&\bar{u}=u+\epsilon\zeta(x,y,t,u)+O(\epsilon^2), \\
& \bar{u}^{(n)}=u^{(n)}+\epsilon\eta^{(n)}+O(\epsilon^2), \quad n\in{\mathbb{N}},  
\end{split}
\end{equation}
for some differentiable functions $\xi, \phi,  \tau, \zeta,$ and $ \eta$ with the corresponding  infinitesimal generator  
\begin{equation}\label{gen}
	X=\xi(x,y,t,u)\frac{\partial}{\partial x}+\phi(x,y,t,u)\frac{\partial}{\partial y}+\tau(x,y,t,u)\frac{\partial}{\partial t}+\zeta(x,y,t,u)\frac{\partial}{\partial u}.
\end{equation}

Thus, the biharmonic heat equation (\ref{key1}) is invariant under the one-parameter Lie group of point transformation (\ref{one}) if and only if it satisfies the invariance criterion
\begin{equation}\label{cri}
X^{[4]}\left[u_t-\left(f^{'}(x)\frac{\partial}{\partial x}+\frac{\partial^2}{\partial x^2}+e^{-2f(x)}\frac{\partial^2}{\partial y^2}\right)^{2}u\right]_{(\ref{key1})}=0.
\end{equation} 
Moreover, the fourth prolongation of the infinitesimal generator (\ref{gen}) is therefore, constructed as follows  
\begin{equation}
	\begin{split}
		X^{[4]}=&X+\eta^{[x]}\frac{\partial}{\partial u_x}+\eta^{[t]}\frac{\partial}{\partial u_t}+\eta^{[xx]}\frac{\partial}{\partial u_{xx}}+\eta^{[xxx]}\frac{\partial}{\partial u_{xxx}}+\eta^{[xxxx]}\frac{\partial}{\partial u_{xxxx}}+\eta^{[yy]}\frac{\partial}{\partial u_{yy}}\\&+\eta^{[xyy]}\frac{\partial}{\partial u_{xyy}}+\eta^{[xxyy]}\frac{\partial}{\partial u_{xxyy}}+\eta^{[yyyy]}\frac{\partial}{\partial u_{yyyy}},
	\end{split}
\end{equation}
which then ensures that the transformation not only applies to the independent and dependent variables $(x,y,t,u)$ but also consistently to all higher-order derivatives in (\ref{key1}); for more details on  prolongation formulae, we refer readers to \cite{ovsiannikov,ibragimov,bluman}.
In addition, upon solving the invariance condition (\ref{cri}), one gets a system of determining equations, which are obtained by equating the coefficients of the various partial derivatives to zero. Hence, the resulting simplified system of determining equations is as follows
\begin{equation*}
\begin{split}
e_1:&  \ \tau_u=\tau_x=\tau_y=\xi_u=\phi_u=\zeta_{uu}=0, \\
e_2:& \ 4\xi_{x}-\tau_{t}=0, \\
e_{3}: &\ \phi_{x}+\xi_{y}{e}^{-2f}=0,  \\
e_{4}: &\ 4\phi_{y}+4\xi\,f_{{x}}-\tau_{{t}}=0, \\
e_{5}: & \ 2\zeta _{{{\it ux}}}-3\xi_{{{xx}}}+\xi f_{{xx}}+f_{{x}}\xi_{{x}}-\xi_{{{yy}}}{{e}^{-2f}}=0,  \\
e_{6}: & \ 2\zeta_{{{uy}}}{{e}^{-2f}}-f_{{x}}\phi_{{x}}+\phi_{{{xx}}}+3\phi_{{{yy}}}{{e}^{-2f}}=0,  \\
e_{7}: & \ -5{e}^{-2f}f_{{x}}\xi_{{y}}+3\phi_{{{xx}}}+{e}^{-2f}\phi_{{{
			yy}}}-2{e}^{-2f}\zeta _{{{uy}}}
+4{e}^{-2f}\xi_{{{xy}}}=0,  \\
e_{8}: & \ \xi_{xx}-2\zeta_{ux}+3{e}^{-2f}\xi_{yy}+3\xi_{x}f_{x}+\left(f_{xx}-2f^2_{x}\right)\xi+4\phi_{xy}-2f_{x}\phi_{y}=0,  \\
e_{9}: & \zeta_{t}-\left(f^{'}(x)\frac{\partial}{\partial x}+\frac{\partial^2}{\partial x^2}+e^{-2f(x)}\frac{\partial^2}{\partial y^2}\right)^{2}\zeta=0, \\
e_{10}:& \ 2f_{x}\zeta_{uy}{{e}^{-2f}}-f_{x}\phi_{yy}{{ e}^{-2f}}+e^{-2f}\left(f^2_{x}-4f_{xx}\right)\xi_{y}+3f_{x}\phi_{xx}-2f_{x}\xi_{xy}{{e}^{-2f}}+2\xi_{xxy}{{e}^{-2f}}-4\zeta_{uxy}{{e}^{-2f}}\\&+2\phi_{xyy}{{e}^{-2f}}+2\phi_{xxx}+2{e}^{-4f}\xi_{yyy}=0, \\
e_{11}:&\ \left(f_{x}f_{xx}-f_{xxx}\right)\xi +2\xi_{xxx}-3f_{x}\zeta _{ux}+3f_{x}\xi_{xx}-\left(f^2_{x}+2
f_{xx}\right)\xi_{x}-f_{x}\xi
_{yy}{{e}^{-2f}}-3\zeta _{uxx}-\zeta _{uyy}{{e}^{-2f}}\\&+2\xi_{xyy}{{e}^{-2f}}=0, \\
e_{12}:&\ 2\phi_{yyy}{{e}^{-2f}}-3\zeta _{uyy}{{e}^{-2f}}-\zeta _{uxx}+2\phi_{xxy}+\left(f_{xxx}+2f^3_{x}-4f_{xx}f_{{x}}\right)\xi+4\left(f_{xx}-f^2_{x}\right)\xi_{x}-2f
_{x}\phi_{xy}+f_{x}\zeta _{ux}\\&+2\left(f^2_{x}-f_{xx}\right)\phi_{y}=0, \\
e_{13}:&\ 2e^{-2f}\left(f^2_{x}-f_{xx}\right)\phi_{yy}-4\zeta _{uxxy}{{e}^{-2f}}+2\phi_{xxyy}{{e}^{-2f}}+4f_{x}\zeta_{uxy}{{ e}^{-2f}}+4e^{-2f}\left(f_{xx}-f^2_{x}\right)\zeta_{uy}-4\zeta_{uyyy}{{e}^{-4f}}\\&+\left(2
f_{xx}+f^2_{x}\right)\phi_{xx}+\phi_{xxxx}+2f_{x}\phi_{xxx}-\phi_{t}-2f_{x}
\phi_{xyy}{{e}^{-2f}}+\left(f_{xxx}+f_{xx}f_{x}\right)\phi_{x}=0, \\
e_{14}:& \ 2\left(2f_{xx}+f^2_{x}\right)\zeta _{ux}+4\zeta_{uxxx}
+\xi_{t}+2e^{-2f}\left(2\zeta_{uxyy}-f_{x}\zeta _{uyy}\right)+6f_{x}\zeta _{uxx}
+4\left(f_{xxx}+f_{xx}f_{x}\right)\xi_{x}\\&+\left(f_
{xxx}f_{x}+ f_{xxxx}+f^2_{xx}\right)\xi-\left(f^{'}(x)\frac{\partial}{\partial x}+\frac{\partial^2}{\partial x^2}+e^{-2f(x)}\frac{\partial^2}{\partial y^2}\right)^{2}\xi=0. 
\end{split}
\end{equation*}
Additionally, one gets from $e_2$ the following 
\begin{equation*}
e_{15}:\ \xi_{xx}=\xi_{xy}=0.
\end{equation*}
Next, differentiating $e_{3}$ with respect to (wrt) $y$ and $e_{4}$ wrt $x$ yields  
\begin{equation*}
	\begin{split}&
		e_{16}:\ e^{-2f}\xi_{yy}+\phi_{xy}=0, \quad \text{and}\quad 
		e_{17}:\ \phi_{xy}+f_x\xi_x+f_{xx}\xi=0, 
	\end{split}
\end{equation*}
while subtracting $e_{16}$ from $e_{17}$ gives 
\begin{equation*}
	e_{18}:\ \xi_{yy}e^{-2f}-\xi_xf_x-\xi f_{xx}=0.
\end{equation*}
Further, substituting $e_{18}$ into $e_{5}$ using $e_2$ yields 
\begin{equation*}
	e_{19}:\ \zeta_{ux}=0. 
\end{equation*}
In the same fashion, differentiating $e_{3}$ wrt $x$ and $e_{4}$ wrt $y$,  one gets 
\begin{equation*}
	\begin{split}&
		e_{20}:\ -2f_xe^{-2f}\xi_{y}+\phi_{xx}=0, \quad \text{and}\quad
		e_{21}:\ \phi_{yy}+f_x\xi_y=0. 
	\end{split}
\end{equation*}
Furthermore, substituting $e_{20}$ and $e_{21}$ into $e_{6}$ using $e_{3}$ reveals 
\begin{equation*}
	e_{22}:\ \zeta_{uy}=0. 
\end{equation*} 
In addition, on differentiating $e_{18}$ wrt to $x,$ one gets 
\begin{equation*}
	e_{23}:\ 2\xi_{yy}e^{-2f}f_x+2\xi_xf_{xx}+\xi f_{xxx}=0, 
\end{equation*}
while eliminate $\xi$ and $\xi_x,$ using $e_{18}$ and $e_{23},$ respectively, gives 
\begin{equation}\label{key5}
	e_{24}: \ e^{-2f}(f_{xxx}+2f_{x}f_{xx})\xi_{yy}-(f_x f_{xxx}-2f^2_{xx})\xi_x=0,
\end{equation}
and 
\begin{equation}\label{key4}
	e_{25}: \ 2e^{-2f}(f_{xx}+f^2_x)\xi_{yy}+(f_{xxx}f_x-2 f^2_{xx})\xi =0. 
\end{equation}
In the same way, eliminating $\xi_{yy}$ using $e_{24}$ and $e_{25}$ reveals 
\begin{equation}\label{key3}
	e_{26}: \ 2(f_{xx}+f^2_x)\xi_x+(f_{xxx}+2f_{x}f_{xx})\xi =0, 
\end{equation}
while differentiating $e_{26}$ with respect to $y$ yields 
\begin{equation}\label{key2}
	e_{27}: \ (f_{xxx}+2f_{x}f_{xx})\xi_y =0. 
\end{equation}

\subsection{$\xi_y=0$}
From \( e_{25} \), one acquires 
\begin{equation}
e_{28}:\ (f_{xxx}f_x-2 f^2_{xx})\xi =0. 
\end{equation}

\subsubsection{$\xi_y=0$ and $\xi=0$} 
The temporal infinitesimal obtained from $e_2$ is
\begin{equation}
	\tau = k_{1}.
\end{equation}
Similarly, from $e_1, e_3, e_4$ and $e_{12}$, one notes ${{\phi }_{u}}={{\phi }_{x}}={{\phi }_{y}}={{\phi }_{t}}=0$, therefore 
\begin{equation}
	\phi = k_{2}.
\end{equation}
Lastly, differentiating $e_{9}$ wrt to $u$, coupled with $e_{1}, e_{19},$ and $e_{22}$, $${{\zeta }_{uu}}={{\zeta }_{ux}}={{\zeta }_{uy}}={{\zeta }_{ut}}=0,$$ and in what follows
\begin{equation}
	\zeta=k_{3}u+\mathcal{G}( x,y,t), 
\end{equation} 
where $\mathcal{G}( x,y,t)$ satisfies $e_9$. This gives the minimal symmetry algebra, for any arbitrary function \( f(x) \) as follows 
\begin{equation}
	X_1 = \frac{\partial}{\partial y}, \quad X_2 = \frac{\partial}{\partial t}, \quad X_3 = u\frac{\partial}{\partial u}, \quad X_u = \mathcal{G} \frac{\partial}{\partial u}.
\end{equation}
Moreover, to get hold of the determining functions \( f(x) \) which may give a larger symmetry algebra, one considers the cases \( \xi_y \neq 0 \).

\subsubsection{$\xi_y=0$ and $\xi\neq{0}$} 
It follows from $e_{28}$ that 
\begin{equation}\label{2.18}
f_{xxx}f_x-2 f^2_{xx}=0. 
\end{equation}
Solving \eqref{2.18} yields 
\begin{equation}
f(x)=\ln \left[ {{\beta }_{5}}{{(x-{{\alpha }_{2}})}^{{{\alpha }_{3}}}} \right], 
\end{equation}
where ${{\beta }_{5}}>0,\;x>{{\alpha }_{2}},\; {{\alpha }_{3}}\ne 0,1$. Additionally, in what follows, 
\[w(x)={{\beta }_{5}}{{(x-{{\alpha }_{2}})}^{{{\alpha }_{3}}}}.\]
It is straightforward to derive from \eqref{2.1} that the Gaussian curvature \(\mathcal{K}\) is then obtained as follows 
\begin{equation}\label{car3}
\mathcal{K} = \frac{\alpha_3 (\alpha_3 - 1)}{(x - \alpha_2)^2}.
\end{equation}
This expression shows that the surface of revolution has a positive variable curvature when \( \alpha_3 < 0 \) or \( \alpha_3 > 1 \). For \( 0 < \alpha_3 < 1 \), the surface exhibits a negative variable curvature. 

\begin{figure}[H]
	\centering
	\begin{minipage}{0.42\textwidth}
		\centering
		\includegraphics[width=\linewidth]{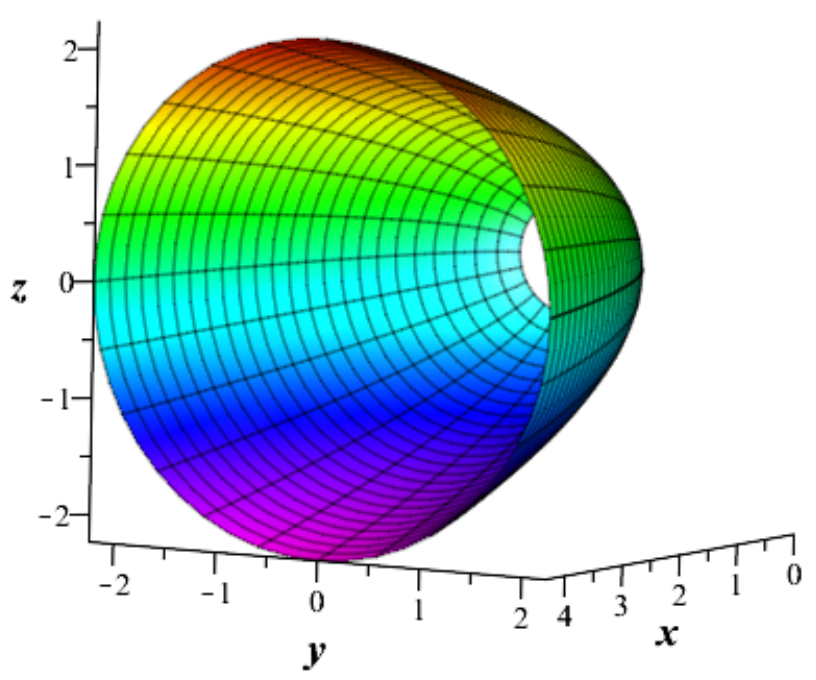} \\
	 \hspace*{-12mm} (a) $\alpha_{2}=0$, $\alpha_{3}=\frac{1}{2},$ $\beta_{5}=1$.
	\end{minipage}
	\hfill
	\begin{minipage}{0.42\textwidth}
		\centering
		\includegraphics[width=\linewidth]{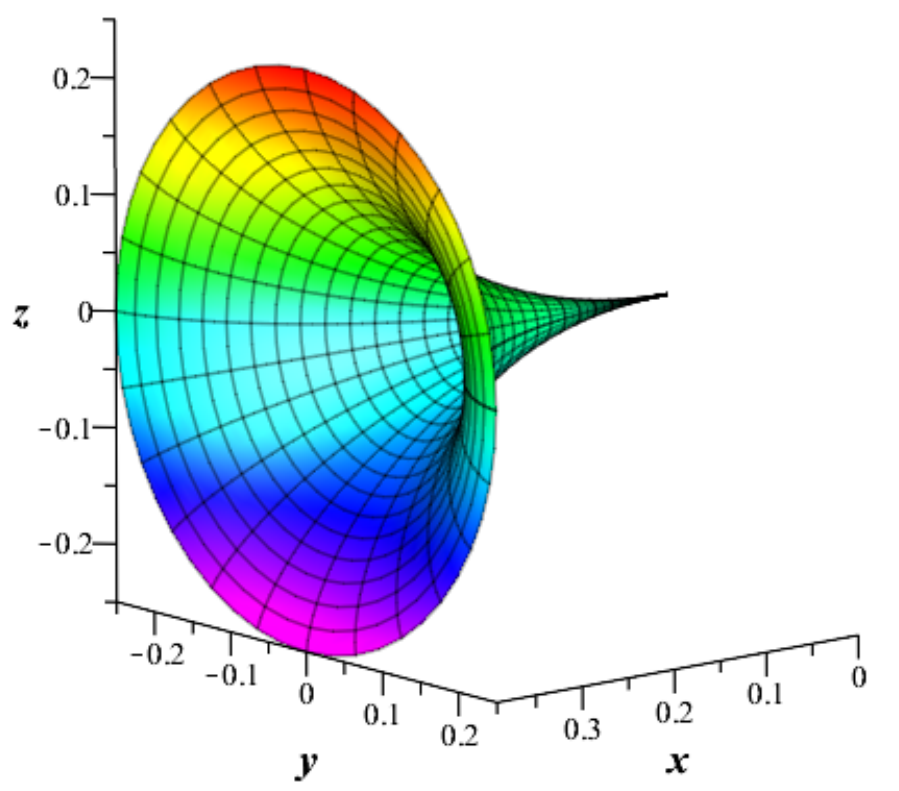}\\
		 \hspace*{-12mm} 
		(b) $\alpha_{2}=0$, $\alpha_{3}={2},$ $\beta_{5}=1.$
	\end{minipage}
	\caption{ Surface of revolution  illustrations when (a) 
	  $\frac{1}{4}\leqslant x\leqslant5$, and (b) 
	$0\leqslant x\leqslant\frac{1}{2}$, for  $0\leqslant y<2\pi.$ }
			\label{fig:image11}
	\end{figure}

Accordingly, solving for $v$ using \eqref{unit}, one thus obtains 
\begin{equation}
v(x)=\int_{{{\alpha }_{4}}}^{x}{\sqrt{1-\beta _{5}^{2}\alpha _{3}^{2}{{(s-{{\alpha }_{2}})}^{2({{\alpha }_{3}}-1)}}}\text{d}s},\quad x\ge {{\alpha }_{4}} , 
\end{equation}
where \[
x \in \left[ \alpha_2 - \left( \frac{1}{\beta_5^2 \alpha_3^2} \right)^{\frac{1}{2\alpha_3 - 2}},\ 
\alpha_2 + \left( \frac{1}{\beta_5^2 \alpha_3^2} \right)^{\frac{1}{2\alpha_3 - 2}} \right]
.\]
Therefore, the parameterization of $\varphi$ is given by 
\begin{equation}
\varphi (x)=\left( \int_{{{\alpha }_{4}}}^{x}{\sqrt{1-\beta _{5}^{2}\alpha _{3}^{2}{{(s-{{\alpha }_{2}})}^{2({{\alpha }_{3}}-1)}}}\text{d}s,}\ \ {{\beta }_{5}}{{(x-{{\alpha }_{2}})}^{{{\alpha }_{3}}}} \right),\quad x\ge {{\alpha }_{4}},
\end{equation}
and thus the corresponding coordinate patch takes the form
\begin{equation}
Y(x,y)=\left( \int_{{{\alpha }_{4}}}^{x}{\sqrt{1-\beta _{5}^{2}\alpha _{3}^{2}{{(s-{{\alpha }_{2}})}^{2({{\alpha }_{3}}-1)}}}\text{d}s},\ \ {{\beta }_{5}}{{(x-{{\alpha }_{2}})}^{{{\alpha }_{3}}}}\cos y,\ \,{{\beta }_{5}}{{(x-{{\alpha }_{2}})}^{{{\alpha }_{3}}}}\sin y \right),
\end{equation}
where $x\ge {{\alpha }_{4}}$ and $0 \leq y < 2\pi;$ see Figure 
		\ref{fig:image11} for the illustrations of such surfaces of revolution.

Furthermore, on solving the determining equations, through the utilization of the determining function $f(x)=\ln \left[ {{\beta }_{5}}{{(x-{{\alpha }_{2}})}^{{{\alpha }_{3}}}} \right], (\alpha_3 \neq{0,1}),$ the symmetry algebra admitted by the biharmonic heat equation on the examining surface is noted to be of finite dimension, which is further  generated by the following vector fields
\begin{equation}
X_1, \quad X_2, \quad X_3, 
\quad X_4=\frac{x-\alpha_{2}}{4}\frac{\partial}{\partial x}-\frac{(\alpha_3-1)}{4}y\frac{\partial}{\partial y}+t\frac{\partial}{\partial t},  \quad X_u.
\end{equation}

\subsection{$\xi_y \neq 0$}
Observe that if $\xi_y\neq 0$, then what follows from \eqref{key2} is that 
\begin{equation}\label{class1}
	f_{xxx}+2f_{x}f_{xx} =0, 
\end{equation}
which implies that 
\begin{equation}\label{2.25}
f_{xx}+f^2_{x}=\eta. 
\end{equation}

\subsubsection{$\eta=0=\mathcal{K}$}
From \eqref{2.25}, one gets 
\begin{equation}\label{class1.12}
f_{xx} + f_x^2 = 0,
\end{equation} 
while \eqref{key4} yields 
\begin{equation}\label{class1.13}
f_{xxx}f_x-2 f^2_{xx} =0.
\end{equation} 
Next, solving for $f$ from the latter equation gives
\begin{equation}
f(x)=\ln|\beta_{3}x+\beta_{4}|,
\end{equation}
where $\beta_4$ and $\beta_4$ are arbitrary constants. What is more, it is worth noting here that the acquired determining function in the latter equation provides all the three flat surfaces of revolution. For an instance, when $ \beta_{3} = 0,$ and $ \beta_{4} = 1$, one respectively obtains a cylinder and a plane; otherwise, the cone surface is revealed \cite{mpungu}. 
\\
\\
\textbf{Cylinder:} Substituting $f(x)=\ln(\beta_4),$ $\beta_4>0$ into the determining equations, the biharmonic heat equation on cylinder admitted the following Lie symmetry algebras
\begin{equation}
\begin{split}
& \qquad \qquad
X_1, \quad X_2, \quad X_3, \quad X_u, 
\quad X_5=\frac{\partial}{\partial x}, \\& X_6=-\beta^{2}_4 y\frac{\partial}{\partial x}+x\frac{\partial}{\partial y}, \quad X_7=\frac{1}{4}x\frac{\partial}{\partial x}+\frac{1}{4}y\frac{\partial}{\partial y}+t\frac{\partial}{\partial t}.
\end{split}
\end{equation}
Thus, a cylindrical surface is, therefore, parameterized as follows 
\begin{equation*}
X(x,y)=\left(x,\beta_4\cos(y),\beta_4\sin(y)\right), \quad 0\le y\le 2\pi
.
\end{equation*}
\textbf{Plane:}
Solving the acquired determining equations for $f(x)=\ln(x+\beta_4),$ $\beta_{4}>0,$ the biharmonic heat equation on a plane admits the following Lie algebra of point symmetries
\begin{equation}
\begin{split}&
\qquad \quad X_1, \quad X_2, \quad X_3, \quad X_u, 
\quad X_8=\frac{x+\beta_4}{4}\frac{\partial}{\partial x}+t\frac{\partial}{\partial t}, \\& X_9=-\cos(y)\frac{\partial}{\partial x}+\frac{\sin(y)}{x+\beta_4}\frac{\partial}{\partial y}, \quad X_{10}=\sin(y)\frac{\partial}{\partial x}+\frac{\cos(y)}{x+\beta_4}\frac{\partial}{\partial y},
\end{split}
\end{equation}
together with the parameterization of the surface as follows 
\begin{equation*}
X(x,y)=\left(x,(x+\beta_4)\cos(y),(x+\beta_4)\sin(y)\right),  \quad 0\le y\le 2\pi, \ x>-\beta_{4}.
\end{equation*}
\textbf{Cone:}
Accordingly, a cone surface is parameterized as follows 
\begin{equation*}
X(x,y)=\left(x,l(x+\beta_4)\cos(y),l(x+\beta_4)\sin(y)\right) \quad 0\le y\le 2\pi, \ x>-\beta_{4},
\end{equation*}
and upon utilizing the determining equations, and the function $f(x)=\ln l(x+\beta_4), \ l>0,$ the symmetry algebra is, therefore, generated by the following vector fields
\begin{equation}
\begin{split}
&
\qquad\qquad\qquad\qquad X_1, \quad X_2, \quad X_3, \quad X_u, 
\quad X_8, \\ 
& X_9=-l\cos(ly)\frac{\partial}{\partial x}+\frac{\sin(ly)}{x+\beta_4}\frac{\partial}{\partial y}, \quad X_{10}=l\sin(ly)\frac{\partial}{\partial x}+\frac{\cos(ly)}{x+\beta_4}\frac{\partial}{\partial y}.
\end{split}
\end{equation}

\subsubsection{$\eta>0,$ $\mathcal{K}<0$}
\begin{figure}[H]
	\centering
	\begin{minipage}{0.45\textwidth}
		\centering
		\includegraphics[width=\linewidth]{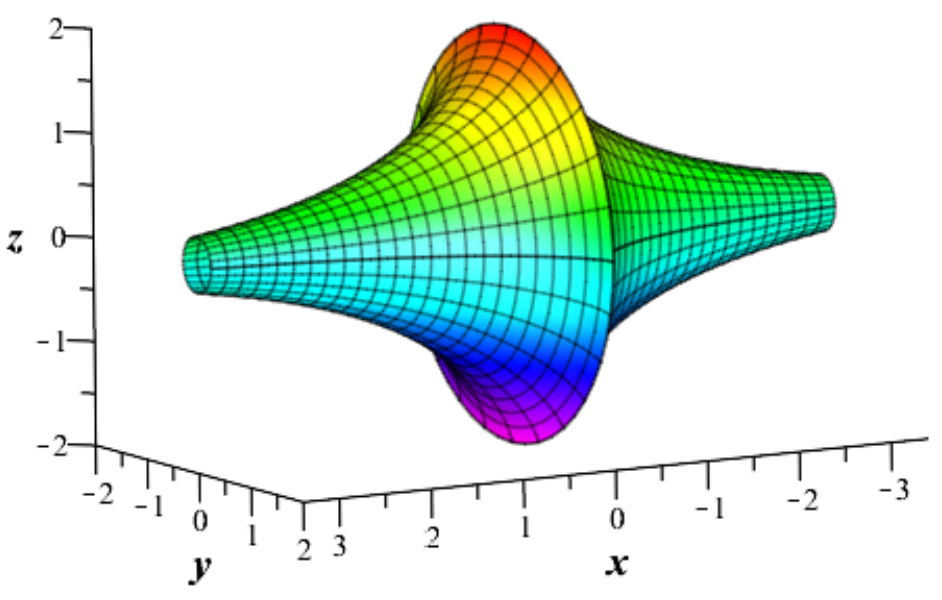} \\
			\end{minipage}
\caption{ A pseudosphere illustration when   
 $\alpha_{2}=0, \alpha_{3}=\frac{1}{2}, \beta_{5}=1, -4\leqslant x\leqslant4, 0\leqslant y<2\pi.$
 }
		\label{fig:image1333}
\end{figure}

\textbf{Case 1: ${{f}_{xx}}=0$ and ${{f}_{x}}\ne 0:$} in this case, solving (\ref{2.25}) yields 
$$f(x) = {\beta_{3}}x + C,$$
for some constants $\beta_{3}\neq 0,$ and $C;$ leading to the expression 
\[w(x)={{\beta }_{4}}{{e}^{{{\beta }_{3}}}}x.\]
Now, setting 
\[{{\alpha }_{3}}=\frac{1}{2{{\beta }_{3}}}\ln \left( \frac{1}{\beta _{4}^{2}\beta _{3}^{2}} \right),\]
it follows from the relationship between \( v(x) \) and \( w(x) \) in \eqref{unit} that 
\[v(x)=\left\{ \begin{split}
	& \int_{{{\alpha }_{3}}}^{x}{\sqrt{1-\beta _{4}^{2}\beta _{3}^{2}{{e}^{2{{\beta }_{3}}s}}}\text{d}s};\quad x\ge {{\alpha }_{3}},\quad \quad {{\beta }_{3}}<0,  \\ 
	& \int_{{{\alpha }_{3}}}^{x}{\sqrt{1-\beta _{4}^{2}\beta _{3}^{2}{{e}^{2{{\beta }_{3}}s}}}\text{d}s};\quad x\le {{\alpha }_{3}},\quad \quad {{\beta }_{3}}>0.   
\end{split} \right.\]
For simplicity, we can assume \( \beta_3 > 0 \), which then allows the function \( \varphi(x) \) to take the form
\[\varphi (x)=\left\{ \begin{split}
	& \left( \int_{{{\alpha }_{3}}}^{x}{\sqrt{1-\beta _{4}^{2}\beta _{3}^{2}{{e}^{-2{{\beta }_{3}}s}}}\text{d}s},\quad {{\beta }_{4}}{{e}^{-{{\beta }_{3}}x}} \right),\quad x\ge {{\alpha }_{3}}, \\ 
	& \left( \int_{{{\alpha }_{3}}}^{x}{\sqrt{1-\beta _{4}^{2}\beta _{3}^{2}{{e}^{2{{\beta }_{3}}s}}}\text{d}s},\quad {{\beta }_{4}}{{e}^{-{{\beta }_{3}}x}} \right),\quad x<{{\alpha }_{3}}.  
\end{split} \right.\]
Certainly, this provides a unit-speed parameterization of the tractrix, and the resulting surface formed by revolving the curve is referred to as a pseudosphere or tractoid \cite{mpungu}. The coordinate representation for this surface is, therefore, expressed as follows 
\[Y(x,y)=\left\{ \begin{split}
	& \left( \int_{{{\alpha }_{3}}}^{x}{\sqrt{1-\beta _{4}^{2}\beta _{3}^{2}{{e}^{-2{{\beta }_{3}}s}}}\text{d}s},\quad {{\beta }_{4}}{{e}^{-{{\beta }_{3}}x}}\cos y,\quad {{\beta }_{4}}{{e}^{-{{\beta }_{3}}x}}\sin y \right),\quad x\ge {{\alpha }_{3}}, \\ 
	& \left( \int_{{{\alpha }_{3}}}^{x}{\sqrt{1-\beta _{4}^{2}\beta _{3}^{2}{{e}^{2{{\beta }_{3}}s}}}\text{d}s},\quad {{\beta }_{4}}{{e}^{{{\beta }_{3}}x}}\cos y,\quad {{\beta }_{4}}{{e}^{{{\beta }_{3}}x}}\sin y \right),\quad x<{{\alpha }_{3}},
\end{split} \right.\]
where $0 \leq y < 2\pi$. Moreover, the Gaussian curvature of the tractoid in this context is a constant negative value $\kappa = -\beta_{3}^2$ \cite{abbena}. In addition, an illustration of a pseudosphere with \( \beta_{4}=2,\;\beta_{3}=\frac{1}{2}\implies \alpha_{3}=0 \), \( -4 \leq x \leq 4 \),  $0\leqslant y<2\pi$ is shown in Figure  \ref{fig:image1333}. 

\begin{figure}[H]
	\centering
	\begin{minipage}{0.45\textwidth}
		\centering
		\includegraphics[width=\linewidth]{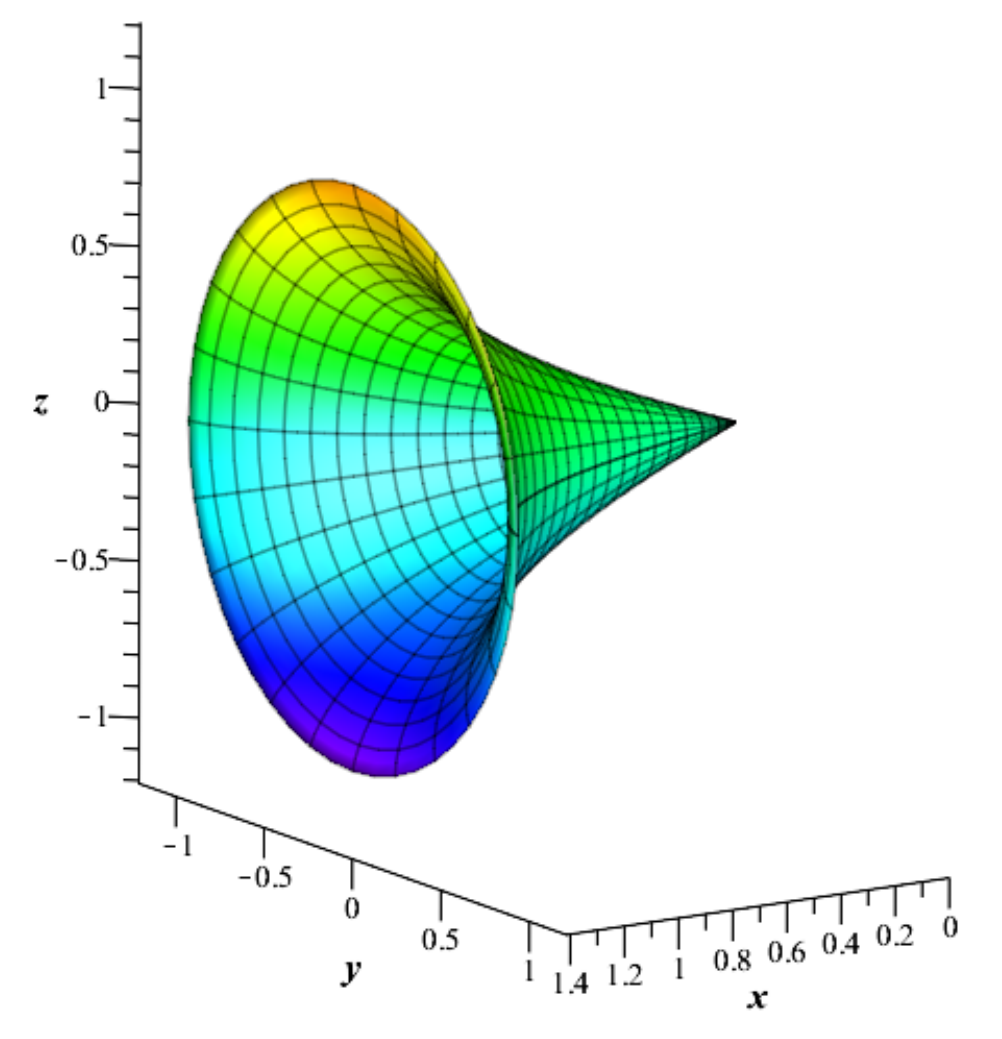}\\
	\end{minipage}
	\caption{
A conic type surface of revolution when 
	$\alpha_{7}=\frac{1}{3},$ $\beta_{7}=1$, $\beta_{8}=0$, $0\leqslant x\leqslant2$, $0\leqslant y<2\pi.$}
			\label{fig:image133}
\end{figure}

Additionally, the symmetry algebra admitted by the corresponding surface is spanned by the following fields 
\begin{equation}
\begin{split}&
X_1, \ \ X_2, \ \  X_3, \   X_u, \
\ X_{11}=-\beta_3^{-1}\frac{\partial}{\partial x}+y\frac{\partial}{\partial y}, \\
& X_{12}=2\beta_3 e^{2C}y\frac{\partial}{\partial x}+\left(e^{-2\beta_{3}x}-\beta^2_{3}e^{2C}y^2\right)\frac{\partial}{\partial y}.
\end{split}
\end{equation}
\textbf{Case 2: ${{f}_{xx}}\neq 0$ and ${{f}_{x}}\ne 0$:}
from \eqref{2.25}, one chooses  $\eta=\beta_{7}^{-2}$, $\beta_{7}>0$, which then leads to 
\begin{equation}\label{class1.17}
f_{xx} + f_x^2 =\beta_{7}^{-2}.
\end{equation}
Setting $g=f_{x}$ in \eqref{class1.17} and solving the resulting equation reveals 
\begin{equation}\label{class1.18}
\ln \left| \frac{{{\beta }_{7}}g+1}{{{\beta }_{7}}g-1} \right|=2\beta _{7}^{-1}(x+c), 
\end{equation}
where $c$ is an arbitrary constant. Next, one firstly considers a case when $g\ge \beta _{7}^{-1}$, in \eqref{class1.18} with ${{\alpha }_{7}}>0,$ yielding 
\begin{equation}\label{class1.19}
f(x)=\ln \left| {{\alpha }_{7}}\sinh \left( \beta _{7}^{-1}x+{{\alpha }_{8}} \right) \right|, 
\end{equation}
where ${{\alpha }_{8}}=c\beta _{7}^{-1}$. In what follows 
\begin{equation}\label{class1.20}
w(x)=\left\{ \begin{matrix}
\ \ {{\alpha }_{7}}\sinh \left( \beta _{7}^{-1}x+{{\alpha }_{8}} \right), \quad x\ge -{{\alpha }_{8}}{{\beta }_{7}},  \\
-{{\alpha }_{7}}\sinh \left( \beta _{7}^{-1}x+{{\alpha }_{8}} \right), \quad x<-{{\alpha }_{8}}{{\beta }_{7}}. 
\end{matrix} \right.
\end{equation}
Consequently, upon considering only the range $ x\ge -{{\alpha }_{8}}{{\beta }_{7}}$, it follows directly from \eqref{unit} that 
\begin{equation}
v(x)=\int_{0}^{\beta _{7}^{-1}x+{{\alpha }_{8}}}{\sqrt{\beta _{7}^{2}-\alpha _{7}^{2}{{\cosh }^{2}}s}}\text{d}s, 
\end{equation}
for $\alpha_{7}\leq \beta_{7},$ and \[-{{\alpha }_{8}}{{\beta }_{7}}<x\leq{{\beta }_{7}}\;\rm{arsinh}\left( \frac{\sqrt{\beta _{7}^{2}-\alpha _{7}^{2}}}{{{\alpha }_{7}}} \right)-{{\alpha }_{8}}{{\beta }_{7}}.\]
The profile surface curve is then deduced from the above relations as follows  
 \[\varphi (x)=\left( \int_{0}^{\beta _{7}^{-1}x+{{\alpha }_{8}}}{\sqrt{\beta _{7}^{2}-\alpha _{7}^{2}{{\cosh }^{2}}s}}\text{d}s,{{\alpha }_{7}}\sinh \left( \beta _{7}^{-1}x+{{\alpha }_{8}} \right) \right),\]
which is then parameterized through the following coordinate patch
\[Y(x,y)=\left( \int_{0}^{\beta _{7}^{-1}x+{{\alpha }_{8}}}{\sqrt{\beta _{7}^{2}-\alpha _{7}^{2}{{\cosh }^{2}}s}}\text{d}s,\ \ {{\alpha }_{7}}\sinh \left( \beta _{7}^{-1}x+{{\alpha }_{8}} \right)\cos y,\ {{\alpha }_{7}}\sinh \left( \beta _{7}^{-1}x+{{\alpha }_{8}} \right)\sin y \right),\]
see Figure \ref{fig:image133} for the graphical illustration of the obtained surface; a conic type surface of revolution curve \cite{abbena}.

\begin{figure}[H]
	\centering
	\begin{minipage}{0.45\textwidth}
		\centering
		\includegraphics[width=\linewidth]{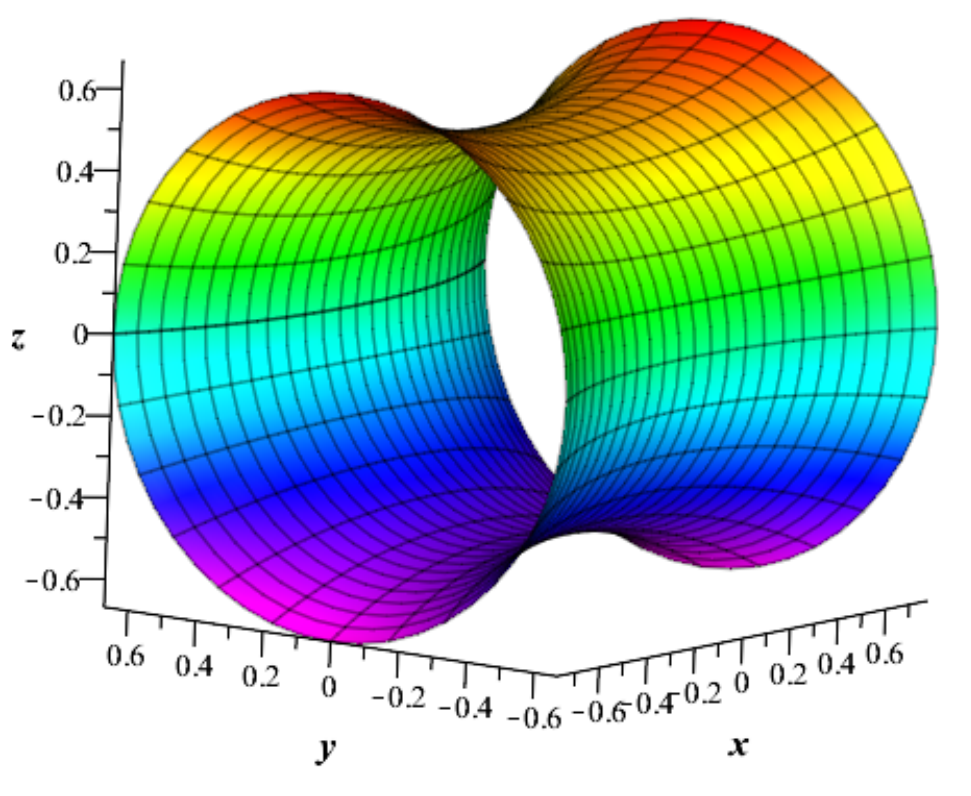}\\
	\end{minipage}
		\caption{
A hyperboloid surface when $\alpha_{7}=\frac{1}{3},$ $\beta_{7}=1$, $\beta_{8}=0$, $0\leqslant x\leqslant2$, $0\leqslant y<2\pi.$}
		\label{fig:image155}	
\end{figure}

\noindent It can be easily shown that the Gaussian curvature of this surface is given by $\mathcal{K}=-\beta_{7}^{-2}$ and the biharmonic heat equation on this surface admits  the following vector fields
\begin{equation}
\begin{split}
& X_1, \ X_2, \ X_{12}=-\alpha_7\cos(\beta_7^{-1}\alpha_7y)\frac{\partial}{\partial x}+\coth(\beta_7^{-1}x)\sin(\beta_7^{-1}\alpha_7y)\frac{\partial}{\partial y}, \\
&  X_3, \ X_u, \ X_{14}=\alpha_7\sin(\beta_7^{-1}\alpha_7y)\frac{\partial}{\partial x}+\coth(\beta_7^{-1}x)\cos(\beta_7^{-1}\alpha_7y)\frac{\partial}{\partial y}.
\end{split}
\end{equation}
Next, considering the case when $g< \beta _{7}^{-1}$ in \eqref{class1.18} with ${{\alpha }_{7}}>0,$ one thus gets  
\begin{equation}\label{class1.21}
\ln \left( \frac{{{\beta }_{7}}g+1}{1-{{\beta }_{7}}g} \right)=2\beta _{7}^{-1}(x+c), 
\end{equation}
where $c$ is an arbitrary constant. In addition, going by the fact that $g=f_{x}$, and further solving \eqref{class1.21}, one then obtains 
\begin{equation}\label{class1.22}
f(x)=\ln \left[ {{\alpha }_{7}}\cosh \left( \beta _{7}^{-1}x+{{\alpha }_{8}} \right) \right],
\end{equation}
where ${{\alpha }_{8}}=c\beta _{7}^{-1}$. Clearly, 
\begin{equation}\label{class1.23}
w(x)={{\alpha }_{7}}\cosh \left( \beta _{7}^{-1}x+{{\alpha }_{8}} \right),.
\end{equation}
where $x\in \mathbb{R}.$ Further, tt follows from \eqref{unit}  that
\begin{equation}\label{class1.20}
v(x)=\int_{0}^{\beta _{7}^{-1}x+{{\alpha }_{8}}}{\sqrt{\beta _{7}^{2}-\alpha _{7}^{2}{{\sinh }^{2}}s}}\text{d}s, 
\end{equation}
for \[-{{\beta }_{7}}\left( \text{arcsinh}\frac{{{\beta }_{7}}}{{{\alpha }_{7}}}+{{\alpha }_{8}} \right)\le x\le {{\beta }_{7}}\left( \text{arcsinh}\frac{{{\beta }_{7}}}{{{\alpha }_{7}}}-{{\alpha }_{8}} \right).\]
Accordingly, one obtains the overall profile curve as follows  
\[\varphi (x)=\left( \int_{0}^{\beta _{7}^{-1}x+{{\alpha }_{8}}}{\sqrt{\beta _{7}^{2}-\alpha _{7}^{2}{{\sinh }^{2}}s}}\text{d}s,{{\alpha }_{7}}\cosh \left( \beta _{7}^{-1}x+{{\alpha }_{8}} \right) \right),\]
which is then parameterized using the following coordinate patch
$$Y(x,y)=\left( \int_{0}^{\beta _{7}^{-1}x+{{\alpha }_{8}}}{\sqrt{\beta _{7}^{2}-\alpha _{7}^{2}{{\sinh }^{2}}s}}\text{d}s,\ \ {{\alpha }_{7}}\cosh \left( \beta _{7}^{-1}x+{{\alpha }_{8}} \right)\cos y,\ {{\alpha }_{7}}\cosh \left( \beta _{7}^{-1}x+{{\alpha }_{8}} \right)\sin y \right).$$
Figure \ref{fig:image155} gives the graphical illustration of the profile curve; indeed, a hyperboloid of one sheet type \cite{abbena}.

In the same veins, the Lie algebra of point symmetries admitted by the biharmonic heat equation on this surface is generated through the following fields 
\begin{equation}
\begin{split}&
X_1, \ \ X_2, \ \  X_{14}=-\alpha_{7}e^{\beta_7^{-1}\alpha_7y}\frac{\partial}{\partial x}+\tanh(\beta_7^{-1}x)e^{\beta_7^{-1}\alpha_7y}\frac{\partial}{\partial y}, \\
& X_3, \ X_u, \ X_{15}=\alpha_{7}e^{-\beta_7^{-1}\alpha_7y}\frac{\partial}{\partial x}+\tanh(\beta_7^{-1}x)e^{-\beta_7^{-1}\alpha_7y}\frac{\partial}{\partial y}.
\end{split}
\end{equation}

\subsubsection{$\eta<0,$ $\mathcal{K}>0,$ ${{f}_{xx}}\neq 0$}
From {\eqref{2.25}}, one then chooses  $\eta=-\beta_{6}^{-2}$, $\beta_{6}>0$,  leading to 
\begin{equation}\label{class1.14}
	f_{xx} + f_x^2 = -\beta_{6}^{-2}.
\end{equation}
Accordingly, solving \eqref{class1.14} with ${{\alpha }_{5}}>0,$ one gets 
\begin{equation}\label{class1.15}
	f(x)=\ln \left| {{\alpha }_{5}}\cos \left( \beta _{6}^{-1}x+{{\alpha }_{6}} \right) \right|,
\end{equation}
where ${{\alpha }_{6}}=c\beta _{6}^{-1}>0$. In what follows, 
\begin{equation}\label{class1.16}
	w(x)=\left| {{\alpha }_{5}}\cos \left( \beta _{6}^{-1}x+{{\alpha }_{6}} \right) \right|,
\end{equation}
such that \eqref{unit} yields 
\begin{equation}
	v(x)=\int_{0}^{\beta _{6}^{-1}x+{{\alpha }_{6}}}{\sqrt{\beta _{6}^{2}-\alpha _{5}^{2}{{\sin }^{2}}s}\text{d}s}.
\end{equation}
Consequently, one deduced from the latter relation, concerning the  parameter $x,$ that 
\[\begin{split}
	&\text{if} \;{{\alpha }_{5}}={{\beta }_{6}}:\; \text{then}\;-{{\beta }_{6}}\left( \frac{\pi }{2}+{{\alpha }_{6}} \right)\le x\le {{\beta }_{6}}\left( \frac{\pi }{2}-{{\alpha }_{6}} \right), \\ 
	& \text{if} \; {{\alpha }_{5}}<{{\beta }_{6}}:\; \text{then}\;x\in \mathbb{R}, \\ 
	& \text{if} \;{{\alpha }_{5}}>{{\beta }_{6}}:\; \text{then}\;-{{\beta }_{6}}\left( {{\alpha }_{6}}+\arcsin \frac{{{\beta }_{6}}}{{{\alpha }_{5}}} \right)\le x\le {{\beta }_{6}}\left( \arcsin \frac{{{\beta }_{6}}}{{{\alpha }_{5}}}-{{\alpha }_{6}} \right). 
\end{split}\]
Therefore, the parameterization of \( \varphi \) can now be expressed as follows
\begin{equation}\label{1.16}
	 \varphi (x)=\left( \int_{0}^{\beta _{6}^{-1}x+{{\alpha }_{6}}}{\sqrt{\beta _{6}^{2}-\alpha _{5}^{2}{{\sin }^{2}}s}\text{d}s},\ \ {{\alpha }_{5}}\cos \left( \beta _{6}^{-1}x+{{\alpha }_{6}} \right) \right), 
\end{equation}
where the corresponding coordinate patch takes the form
\begin{equation}\label{1.17}
Y(x,y)=\left( \int_{0}^{\beta _{6}^{-1}x+{{\alpha }_{6}}}{\sqrt{\beta _{6}^{2}-\alpha _{5}^{2}{{\sin }^{2}}s}\text{d}s},\ \ {{\alpha }_{5}}\cos \left( \beta _{6}^{-1}x+{{\alpha }_{6}} \right)\cos y,\ {{\alpha }_{5}}\cos \left( \beta _{6}^{-1}x+{{\alpha }_{6}} \right)\sin y \right), 
\end{equation}
for $0 \leq y < 2\pi$. Notably, the acquired surfaces in the latter equation possessed positive Gaussian curvatures, that is,  $\mathcal{K}=\beta_{6}^2.$ In addition, the symmetry algebra admitted by the governing heat equation on the examining surface is thus spanned by the following vector fields
\begin{equation}
\begin{split}
&
X_1, \ \ \ X_2, 
\quad X_{16}=\alpha_5\cos(\beta_6^{-1}\alpha_5y)\frac{\partial}{\partial x}+\tan(\beta_6^{-1}x)\sin(\beta_6^{-1}\alpha_5y)\frac{\partial}{\partial y}, \\
& 
X_3, \quad X_u,  \ X_{17}=-\alpha_5\sin(\beta_6^{-1}\alpha_5y)\frac{\partial}{\partial x}+\tan(\beta_6^{-1}x)\cos(\beta_6^{-1}\alpha_5y)\frac{\partial}{\partial y}.
\end{split}
\end{equation}

\begin{remark}
Let $ \mathcal{S}(\alpha_{5}, \beta_{6}) $ be the surface of revolution generated by the parameterized unit speed curve as follows \(\varphi(x) = (v, w) \). Then, one obtains a 
	\begin{itemize}
		\item \textbf{Sphere type}: if $\alpha_{5}= \beta_{6}, $ then \( \mathcal{S}(\alpha_{5}, \beta_{6}) \) represents a standard sphere of radius \( \beta_{6} \) as illustrated in Figure \ref{fig:image1rr5}.		
		\begin{figure}[H] 
			\centering
			\begin{minipage}{0.45\textwidth}
				\centering
				\includegraphics[width=\linewidth]{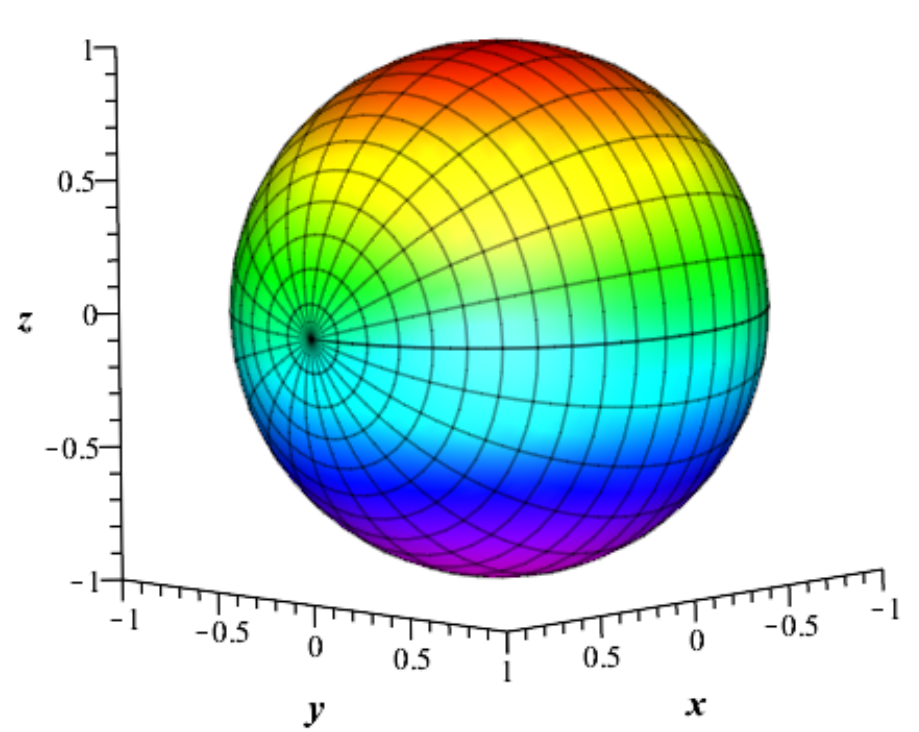}
			\end{minipage}
	\caption{ A sphere type surface when    
	$\alpha_{5}=\beta_{6}=1$, $\beta_{6}=0$, $-\frac{\pi}{2}\leqslant x\leqslant\frac{\pi}{2}$, $0\leqslant y<2\pi$. }
							\label{fig:image1rr5}
								\end{figure}
								
		\item \textbf{Spindle type}: if \( 0 < \alpha_{5}< \beta_{6} \), then \( \mathcal{S}(\alpha_{5}, \beta_{6}) \) is a surface of revolution resembling a rugby ball, with sharp vertices along the axis of revolution; see Figure \ref{fig:image1ff} for the graphical  illustration.
		\begin{figure}[H]
			\centering
			\begin{minipage}{0.45\textwidth}
				\centering
				\includegraphics[width=\linewidth]{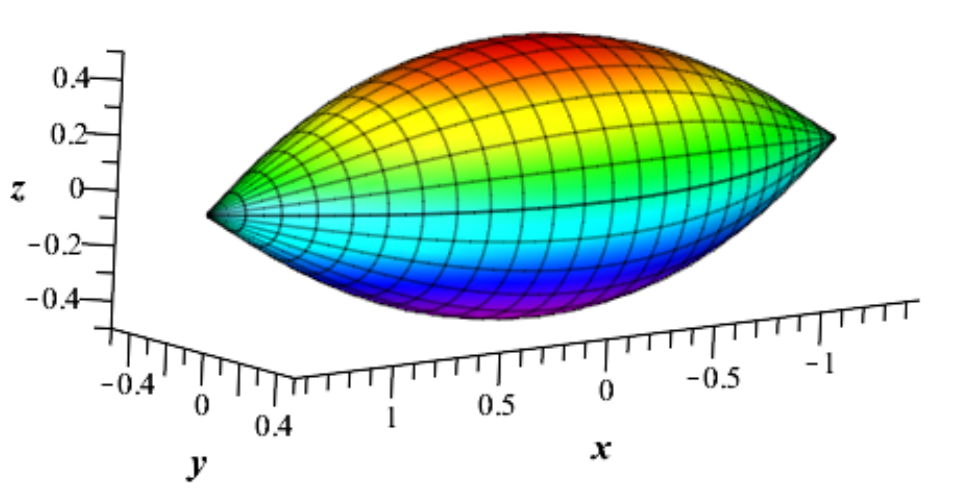}
							\end{minipage}
	\caption{A spindle type surface when $\alpha_{5}=\frac{1}{2}, \beta_{6}=1$, $\beta_{6}=0$, $-1.56\leqslant x\leqslant 1.56$, $0\leqslant y<2\pi.$}
					\label{fig:image1ff}
		\end{figure}
		\item \textbf{Bulge type}: if \( 0 <  \beta_{6}<\alpha_{5} \), then \( \mathcal{S}(\alpha_{5}, \beta_{6})\) is a barrel-shaped surface, which does not intersect the axis of revolution \cite{abbena}; see Figure \ref{fig:image1ffff} for its graphical illustration. 
		\begin{figure}[H]
			\centering
			\begin{minipage}{0.45\textwidth}
				\centering
				\includegraphics[width=\linewidth]{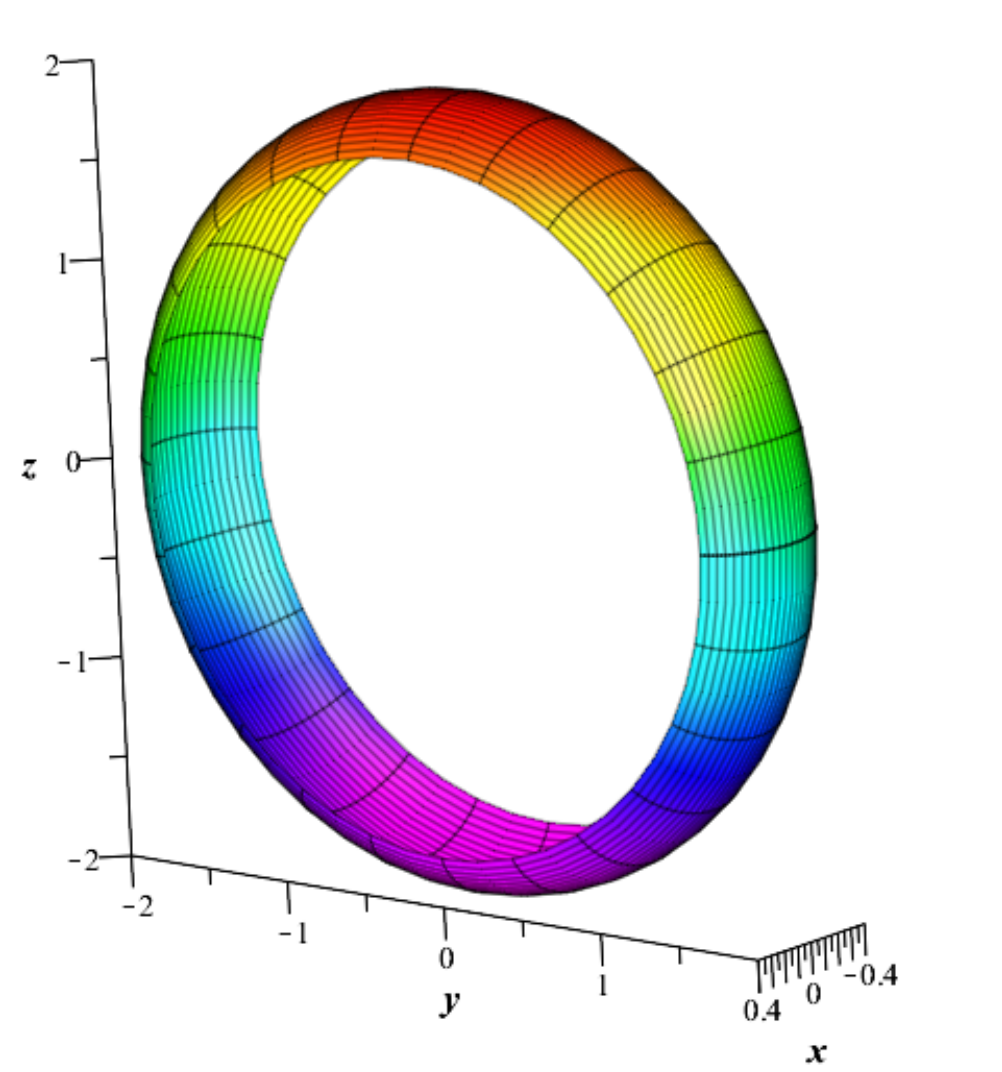}
						\end{minipage}
	\caption{A bulge type surface when $\alpha_{5}=2, \beta_{6}=1$, $\beta_{6}=0$, $-0.5\leqslant x\leqslant 0.5$, $0\leqslant y<2\pi.$}
										\label{fig:image1ffff}
		\end{figure}
	\end{itemize}
\end{remark}

\section{Symmetry Reductions and Exact Solutions}
In this section, the Lie point symmetries admitted by the governing biharmonic heat equation are used to reduce the equation on a surface of revolution. The symmetry reductions will be carried out using the conventional method of introducing similarity variables.

\subsection{Reduction by 1-dimensional subalgebra}

\subsubsection{Subalgebra $\mathcal{Y} = \langle X_1 + aX_2 \rangle, \, a \neq 0$}
  With the consideration of the following subalgebra 
$$X = \frac{\partial}{\partial y} + a\frac{\partial}{\partial t},$$
one thus obtains the resulting characteristic equation as follows 
\[\frac{\text{d}t}{a}=\frac{\text{d}x}{0}=\frac{\text{d}y}{1}=\frac{\text{d}u}{0}, \]
which when solved gives
\begin{equation}\label{2.5}
\eta =x,\text{  }\upsilon =y-\frac{t}{a},\;\text{  and }\;\phi (\eta ,\upsilon )=u.
\end{equation}
Therefore, substituting the latter expressions into \eqref{key1}, one gets
\begin{equation}\label{4.2}
	\begin{aligned}
 -\frac{1}{a}{{\phi }_{\upsilon }}=	&\;{{\phi }_{\eta \eta \eta \eta }}+\left( {f}'''+{f}'{f}'' \right){{\phi }_{\eta }}+\left( 2{f}''+{{{{f}'}}^{2}} \right){{\phi }_{\eta \eta }}+2{f}'{{\phi }_{\eta \eta \eta }} \\ 
	& -2{{e}^{-2f(\eta)}}\left( {f}''-{{{{f}'}}^{2}} \right){{\phi }_{\upsilon \upsilon }}+{{e}^{-4f(\eta)}}{{\phi }_{\upsilon \upsilon \upsilon \upsilon }}. 
\end{aligned}
\end{equation}
Accordingly, substituting the following ansatz   
\begin{equation}
\phi(\eta,\upsilon)=c_1\upsilon+\psi(\eta), 
\end{equation}
into (\ref{4.2}), one thus gets the following reduced equation 
\begin{equation}\label{redu}
\begin{aligned}
&\;{{\psi }_{\eta \eta \eta \eta }}+\left( {f}'''+{f}'{f}'' \right){{\psi }_{\eta }}+\left( 2{f}''+{{{{f}'}}^{2}} \right){{\psi }_{\eta \eta }}+2{f}'{{\psi }_{\eta \eta \eta }}=-\frac{c_1}{a}. 
\end{aligned}
\end{equation}

\subsubsection{Subalgebra  $\mathcal{X}=\langle a{{X}_{2}}+b{{X}_{3}}\rangle \text{  },\text{ }a\ne 0,\text{ }b\ne 0$}
The present consideration thus examines the following subalgebra 
$$X=a\frac{\partial }{\partial t}+bu\frac{\partial }{\partial u},$$
which yields results in the following characteristic equation 
\[ \frac{\text{d}t}{a}=\frac{\text{d}x}{0}=\frac{\text{d}y}{0}=\frac{\text{d}u}{bu}.\]
Consequently, the solution of the above equations further reveals 
\[\eta =x,\;\upsilon =y\;\;\text{and}\;\;\phi(\eta ,\upsilon )=\ln u-\frac{b}{a}t, \]
which then yields the following invariant solution is
\begin{equation}\label{iv}
u(\eta, \upsilon, t) = \exp \left( \phi(\eta, \upsilon) + \frac{b}{a} t \right),
\end{equation}
where \(\eta = x \) and \( \upsilon=y \). Further, upon substituting the invariant solution (\ref{iv}) into (\ref{key1}), one eventually obtains the following transformed equation 
 \begin{equation}
\begin{aligned}
	\frac{b}{a}=\;& 2{{\phi }_{\eta }}{{\phi }_{\eta \eta \eta }}+2\phi _{\eta \eta }^{2}+{{\phi }_{\eta \eta \eta \eta }}+\left( {f}'''+{f}'{f}'' \right){{\phi }_{\eta }}+\left( 2{f}''+{{{{f}'}}^{2}} \right)\left( \phi _{\eta }^{2}+{{\phi }_{\eta \eta }} \right)+2{f}'\left( 2{{\phi }_{\eta }}{{\phi }_{\eta \eta }}+{{\phi }_{\eta \eta \eta }} \right) \\ 
	& -2{{e}^{-2f(\eta )}}\left( {f}''-{{{{f}'}}^{2}} \right)\left( \phi _{\upsilon }^{2}+{{\phi }_{\upsilon \upsilon }} \right)-2{f}'{{e}^{-2f(\eta )}}\left( 2{{\phi }_{\eta }}{{\phi }_{\upsilon }}{{\phi }_{\upsilon \upsilon }}+{{\phi }_{\upsilon \upsilon \upsilon }} \right) \\ 
	& +2{{e}^{-2f(\eta )}}\left( 2{{\phi }_{\eta \eta }}{{\phi }_{\upsilon \upsilon }}+2\phi _{\eta }^{2}{{\phi }_{\upsilon \upsilon }}+2{{\phi }_{\eta }}{{\phi }_{\upsilon \upsilon \upsilon }} \right)+{{e}^{-4f(\eta )}}\left( 2{{\phi }_{\upsilon }}{{\phi }_{\upsilon \upsilon \upsilon }}+{{\phi }_{\upsilon \upsilon \upsilon \upsilon }} \right). \\ 
\end{aligned}
\end{equation}
\subsection{Reduction by 2-dimensional subalgebra}

\subsubsection{Subalgebra $\mathcal{L}=\langle {{X}_{1}},a{{X}_{2}}+b{{X}_{3}}\rangle $}
With the adoptation of the current subalgebra, the resultant characteristic system  corresponding to the subalgebra $${X}_{1}=\frac{\partial }{\partial y},$$ 
is given by
\[\frac{\text{d}t}{0}=\frac{\text{d}x}{0}=\frac{\text{d}y}{1}=\frac{\text{d}u}{0},\]
which then leads to similarity variables
\[\eta =x,\;\upsilon =t,\;\text{and}\;\phi (\eta ,\upsilon )=u.\]
Accordingly, substituting the above similarity variables into \eqref{key1} yields the following reduced equation 
\begin{equation}\label{4.3}
	{{\phi }_{\upsilon }}={{\phi }_{\eta \eta \eta \eta }}+\left( {f}'''+{f}'{f}'' \right){{\phi }_{\eta }}+\text{ }\left( 2{f}''+{{{{f}'}}^{2}} \right){{\phi }_{\eta \eta }}+2{f}'{{\phi }_{\eta \eta \eta }}.
\end{equation}
Notably, since $[{{X}_{1}},a{{X}_{2}}+b{{X}_{3}}]=0$, clearly, the two symmetries commute. Consequently, the second symmetry is inherited by \eqref{4.3}, as it commutes with the first symmetry \cite[p.~285]{ibragimov}, and therefore,
\[Y=a\frac{\partial }{\partial \upsilon }+b\phi \frac{\partial }{\partial \phi },\]
together with the corresponding characteristic system as follows
\[\frac{\text{d}\eta }{0}=\frac{\text{d}\upsilon }{a}=\frac{\text{d}\phi }{b\phi }.\]
In the same way, the latter system is then solved to get hold of the following similarity variables 
 \[r(\eta ,\upsilon )=\eta, \text{ and }w(r)=\ln \phi -\frac{b}{a}\upsilon.\]
Moreover, given the function \( \phi = \exp\left( w(r) + \frac{b}{a}\upsilon \right) \) and \( r(\eta, \upsilon) = \eta \), equation  \eqref{4.3} then gives the following reduced equation  
\begin{equation}
	\begin{aligned}
	 \frac{b}{a}=&\;{{w}^{(4)}}+4{w}'{w}'''+3{{{{w}'}}^{2}}{w}''+{{{{w}'}}^{4}}+\left( {f}'''+{f}'{f}'' \right){w}'+\left( 2{f}''+{{{{f}'}}^{2}} \right)\left( {w}''+{{{{w}'}}^{2}} \right) \\ 
	& +2{f}'\left( {w}'''+3{w}'{w}''+{{{{w}'}}^{3}} \right). 
\end{aligned}
\end{equation}

\subsection{Exact solutions on surfaces of specific curvature}
Here, having obtained all the possible reduced invariant differential equations that emanate form the governing biharmonic heat equation in the above subsections, this subsection then specifically constructs various analytical solutions for the reduced invariant equation (\ref{redu}), considering some specific geometries of interest. Moreover, the chosen reduced differential equation is handier in comparison with the rest of the realized deduced differential equations.

 \begin{examples}[Surface of revolution with zero Gaussian curvature (cylinder)]
 Consider the case where the surface of revolution is a cylinder, corresponding to $f(x)=\ln(\beta).$ Thus, the governing biharmonic heat equation \eqref{key1} on a cylinder admits an exact analytical solution 
 through (\ref{redu}) as follows 
 \begin{equation}
u(x,y,t)=c_1\left(y-\frac{t}{a}-\frac{1}{24a}x^4\right)+\frac{c_2}{6}x^3+\frac{c_3}{2}x^2+c_4x+c_5,
\end{equation}
where $c_j(j=1,2,3,...,5)$ are arbitrary constants.  \\  
\end{examples}

\begin{examples}[Surface of revolution with a constant negative Gaussian curvature (pseudosphere or tractoid)]
Considering the situation where a constant negative Gaussian curvature arisen, which then corresponds to $f(x)=x,$ the governing equation is then said to have a pseudosphere as the involving surface of revolution. Thus, the governing heat model in \eqref{key1} through (\ref{redu}) satisfies the following analytical solution 
 	\begin{equation}
u(x,y,t)=c_1\left(y-\frac{t}{a}-\frac{1}{2a}x^2\right)+(c_2(x+2)+c_3)e^{-x}+c_4x+c_5,  
	\end{equation}    
where $c_j(j=1,2,3,...,5)$ are arbitrary constants.  \\  

\end{examples}

\begin{examples}[Surface of revolution with a constant positive Gaussian curvature (paraboloid)]
	For a paraboloid surface of revolution, the corresponding function $f(x)$ takes the following form $f(x)=\frac{1}{2}\ln(x), \ x>0.$ In addition,  the governing heat model in \eqref{key1} through (\ref{redu}) when a paraboloid surface is considered satisfies the following analytical solution 
	\begin{equation}
u(x,y,t)=c_1\left(y-\frac{t}{a}-\frac{1}{42a}x^4\right)+\frac{1}{2}c_2x^2+2c_3\sqrt{x}+\frac{2}{5}c_4x^{\frac{5}{2}}+c_5, 
	\end{equation}
where $c_j(j=1,2,3,...,5)$ are arbitrary constants.  \\

\end{examples}

\section{Conclusion}
In this paper, we conducted a complete Lie symmetry group classification of the biharmonic heat equation on a generalized surface of revolution. We classified the infinitesimal symmetries associated with the equation and demonstrated that the symmetry algebra remained structurally identical to that of the harmonic heat equation on such surfaces \cite{mpungu}. The result of the study further revealed profound geometric and analytical connections between the harmonic and biharmonic heat equations.  In addition, the study further utilized the constructed symmetries to obtain similarity reductions, transforming the original partial differential equation into a lower-order and more manageable class of differential equations. Exact solutions were derived for surfaces with zero, positive, and negative Gaussian curvature. Lastly, the overall analysis shows the power of Lie symmetry methods in simplifying and solving complex differential equations on curved geometries.

\bibliographystyle{abbrv}
\bibliography{refbih}
\end{document}